\documentclass[12pt]{amsart}
\hoffset= - 0.75 in

\pdfoutput=1

\usepackage{amssymb,amsmath,amsfonts,epsfig,latexsym}

\newtheorem*{Richman}{Richman's Theorem}

\theoremstyle{definition}

\def\<{\ensuremath{\langle}}
\def\>{\ensuremath{\rangle}}

\begin{document}

\title{Bidding chess}

\author[Bhat]{Jay Bhat}
\author[Payne]{Sam Payne}

\maketitle

It all started with a chessboard
\vspace{-30pt}
\begin{center}
\scalebox{.3}{\includegraphics{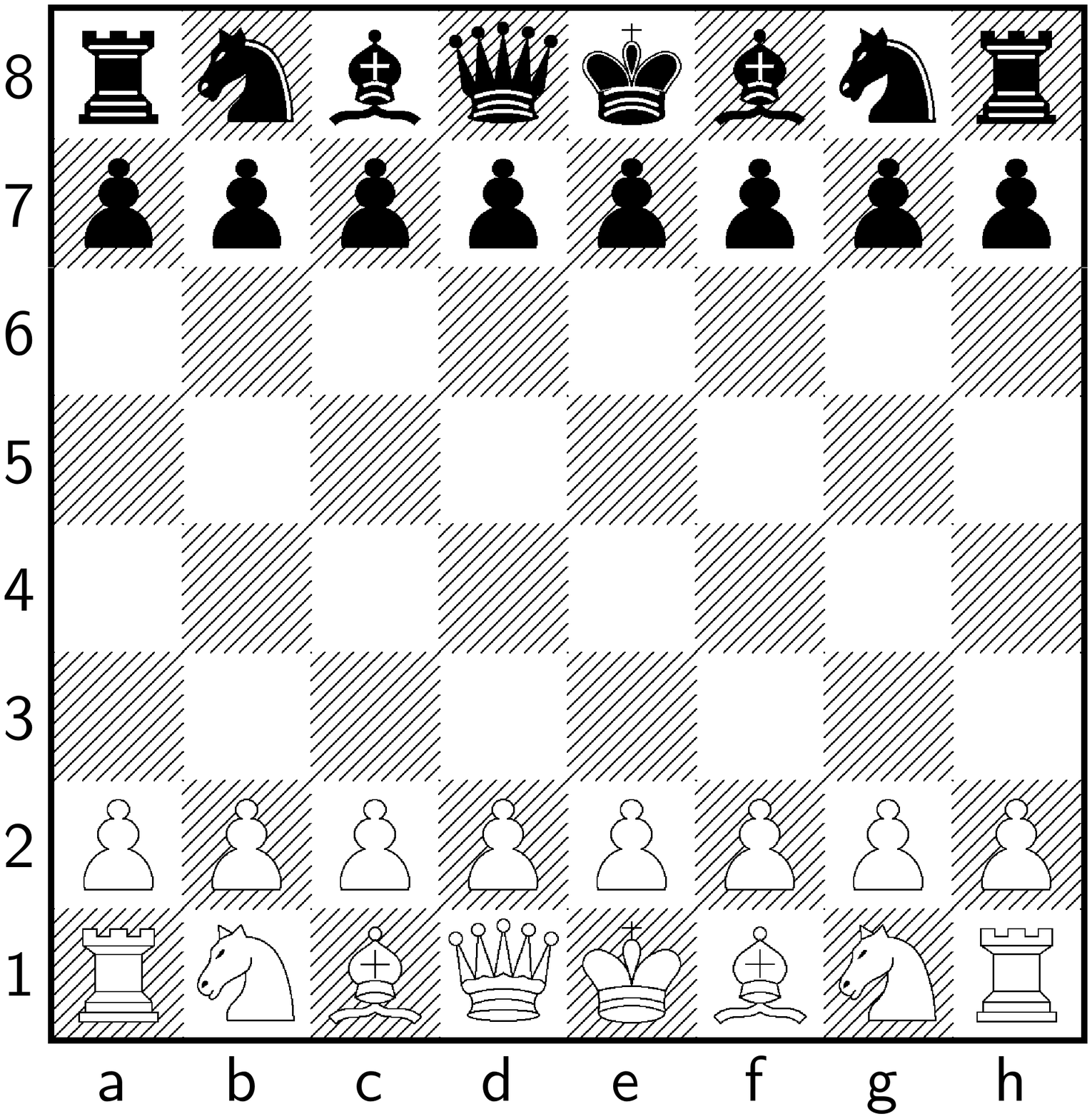}}
\end{center}
\vspace{-30pt}
\noindent and a bottle of raki

\begin{center}
\scalebox{.4}{\includegraphics{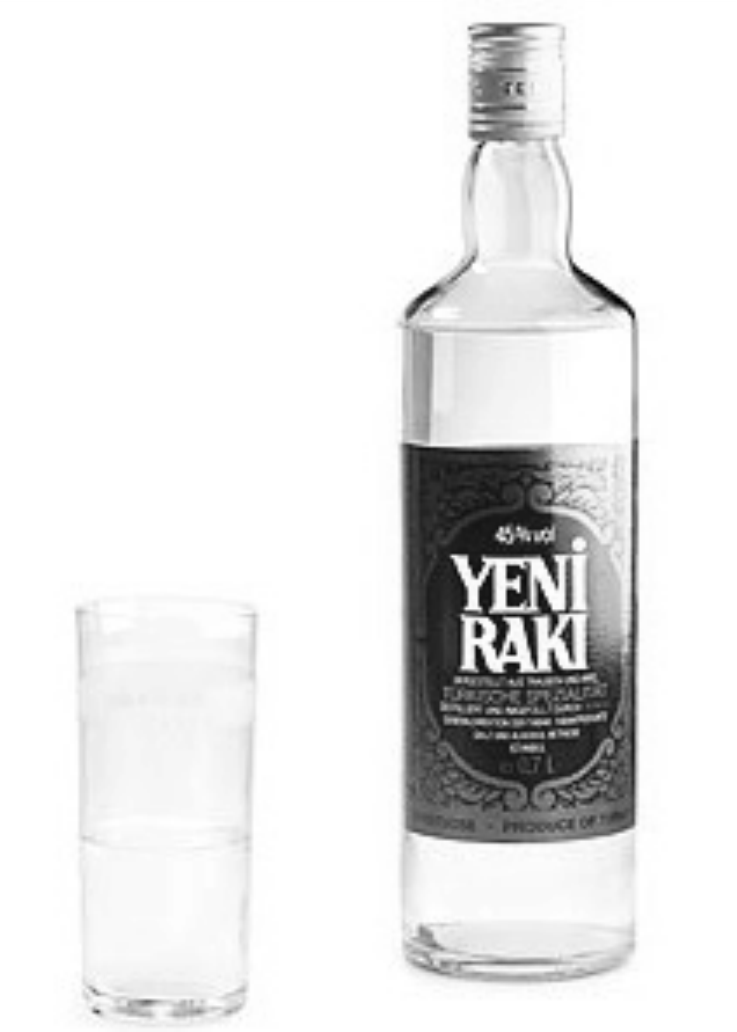}}
\end{center}

\noindent at an otherwise respectable educational institution.  SP's friend Ed had just returned from Turkey with a bottle of the good stuff.  We were two and a half sheets to the wind, feeling oppressed by the freshman dormroom, its cinderblock walls, and our mediocre chess skills.  The details of what transpired are clouded with alcohol and time, but Ed was enrolled in a seminar on auction theory and a site called eBay had recently emerged on the web.  Somehow we decided that alternating moves was boring---too predictable.  Some moves are clearly worth more than others, and the game would be much more interesting if you had to pay for the right to move.  Enter Bidding Chess.

{\bf How to play.}  We both start with one hundred chips.  Before each move, we write down our bids, and the player who bids more gives that many chips to the other player and makes a move on the chessboard.  For example, if I bid nineteen for the first move and you bid twenty-four, then you give me twenty-four chips and make a move on the chessboard.  Now I have 124 chips and you have 76 and we bid for the second move.  By the way, you bid way too much and now you're toast!

{\bf How to win.} You win by capturing your opponent's king.  Or, rather, I win by capturing your king.  None of this woo-woo checkmate stuff.  I don't care if you have me in checkmate when I have enough chips to make seven moves in a row.

{\bf What if the bids are tied?} Quibblers.  The bids are never tied.  There is an extra chip, called *.  If you have * in your pile, then you can include it with your bid and win any tie.  And if you win, then * goes to your opponent along with the rest of your bid.  But if you wuss out and save * for later, then you lose the tie.

Bidding Chess is meant to be played, so set up the board and grab a friend, maybe one you never really liked much anyway, and try it out.  Think carefully----the game is never won on the first move, but it is often lost there.  After you have played a few times, take a look at the following transcript from a game played in 1015 Evans Hall, the common room of the UC Berkeley math department, in fall 2006.  Names have been changed for reasons you may imagine.  We write $N$* to denote a pile or bid with $N$ chips plus the * chip.

{\bf A sample game.}  Alice and Bob both start with one hundred chips.  Alice offers Bob the * chip, but he refuses.  Alice shrugs, takes the black pieces,  and bids twelve for the first move.  Bob bids thirteen and moves his knight to c3.  
\vspace{-30pt}
\begin{center}
\scalebox{.3}{\includegraphics{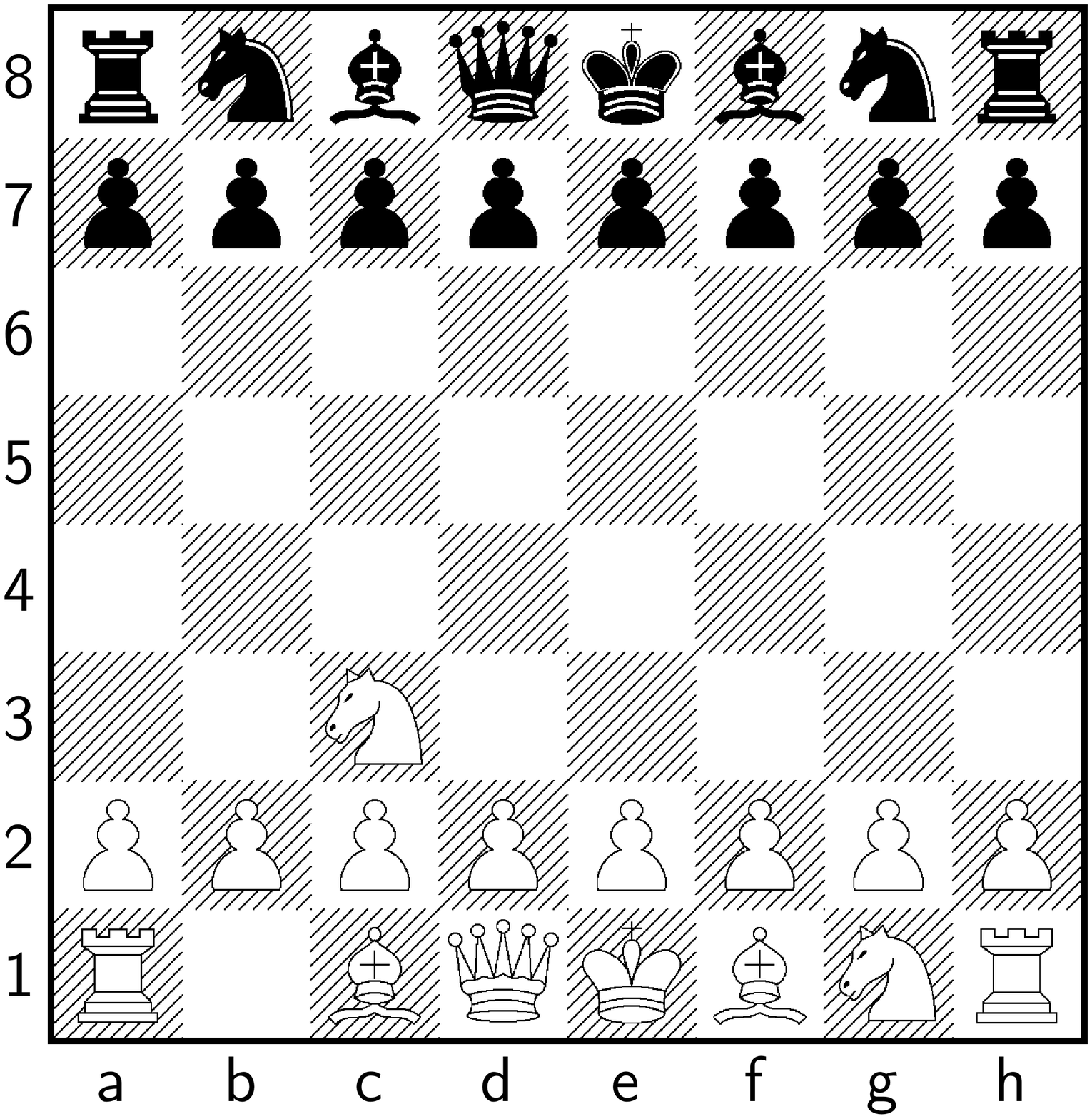}}
\end{center}
\vspace{-30pt}
Now Alice has 113* chips, and Bob has 87 as they ponder the value of the next move.  Alice figures that the second move cannot possibly be worth more than the first, because it would be silly to bid more than thirteen and end up in a symmetric position with fewer chips than Bob.  So she bids eleven plus *, which feels about right.  Bob reasons similarly, and also bids eleven.  Alice wins the tie with * and moves her king's pawn forward one to e6.

Bob, who played competitive chess in high school, is puzzled by this conservative opening move.  Feeling comfortable with his board position, he decides to bid only nine of his 98* chips for the third move and is mildly surprised when Alice bids fifteen.  She moves her bishop to c5.
\vspace{-30pt}
\begin{center}
\scalebox{.3}{\includegraphics{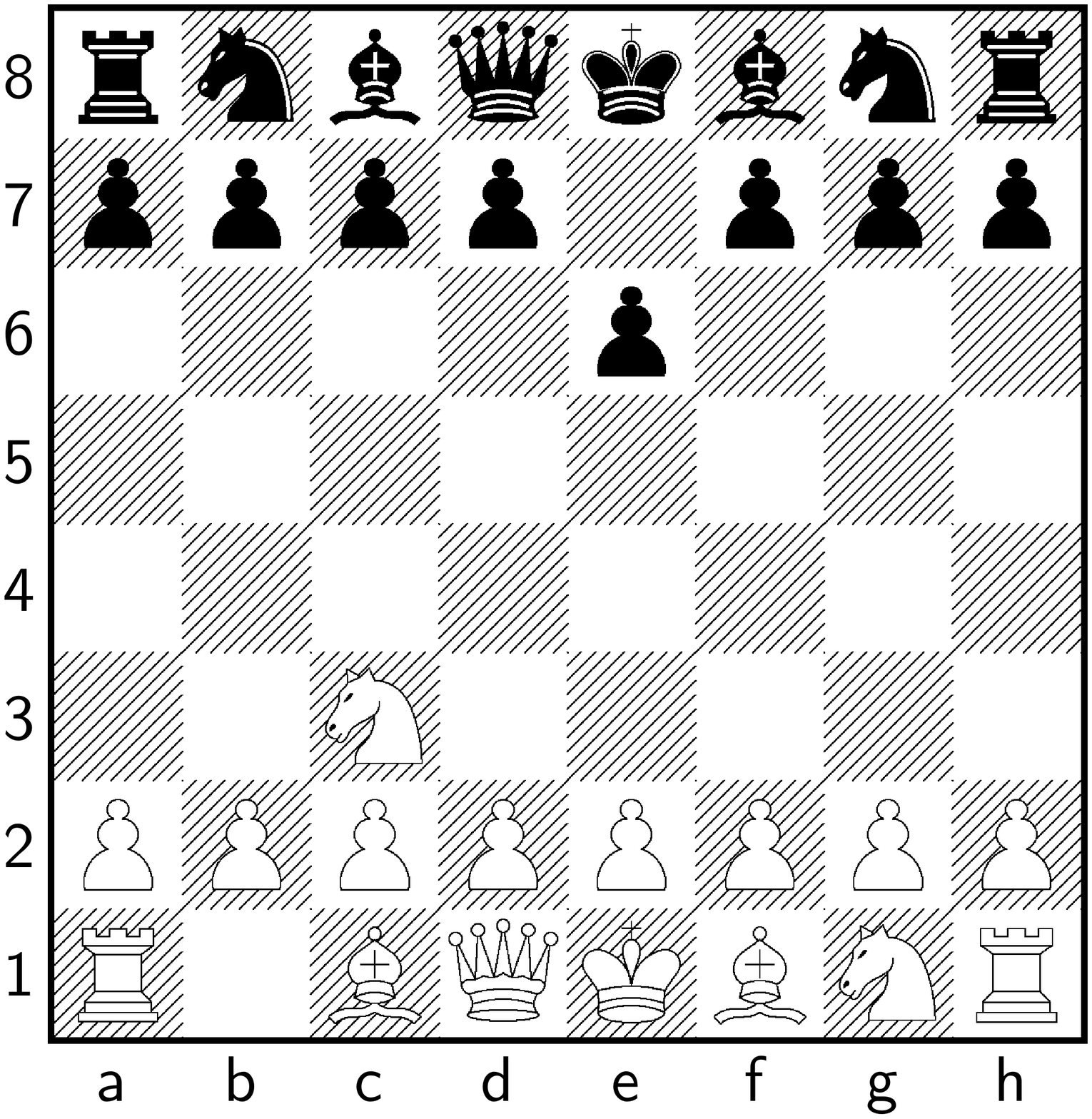}}
\end{center}
\vspace{-30pt}
Now Alice has 87 chips, while Bob has 113*.  Since Alice bid 15 for the last move and started an attack that he would like to counter, Bob bids fifteen for the next move, which seems fair.  Alice bids twenty-two and takes the pawn at f2.

 Bob realizes with some dismay that he must win the next move to prevent Alice from taking his kind.  She bids all 65 of her chips, and he bids 65*.  King takes bishop.
\vspace{-30pt}
\begin{center}
\scalebox{.3}{\includegraphics{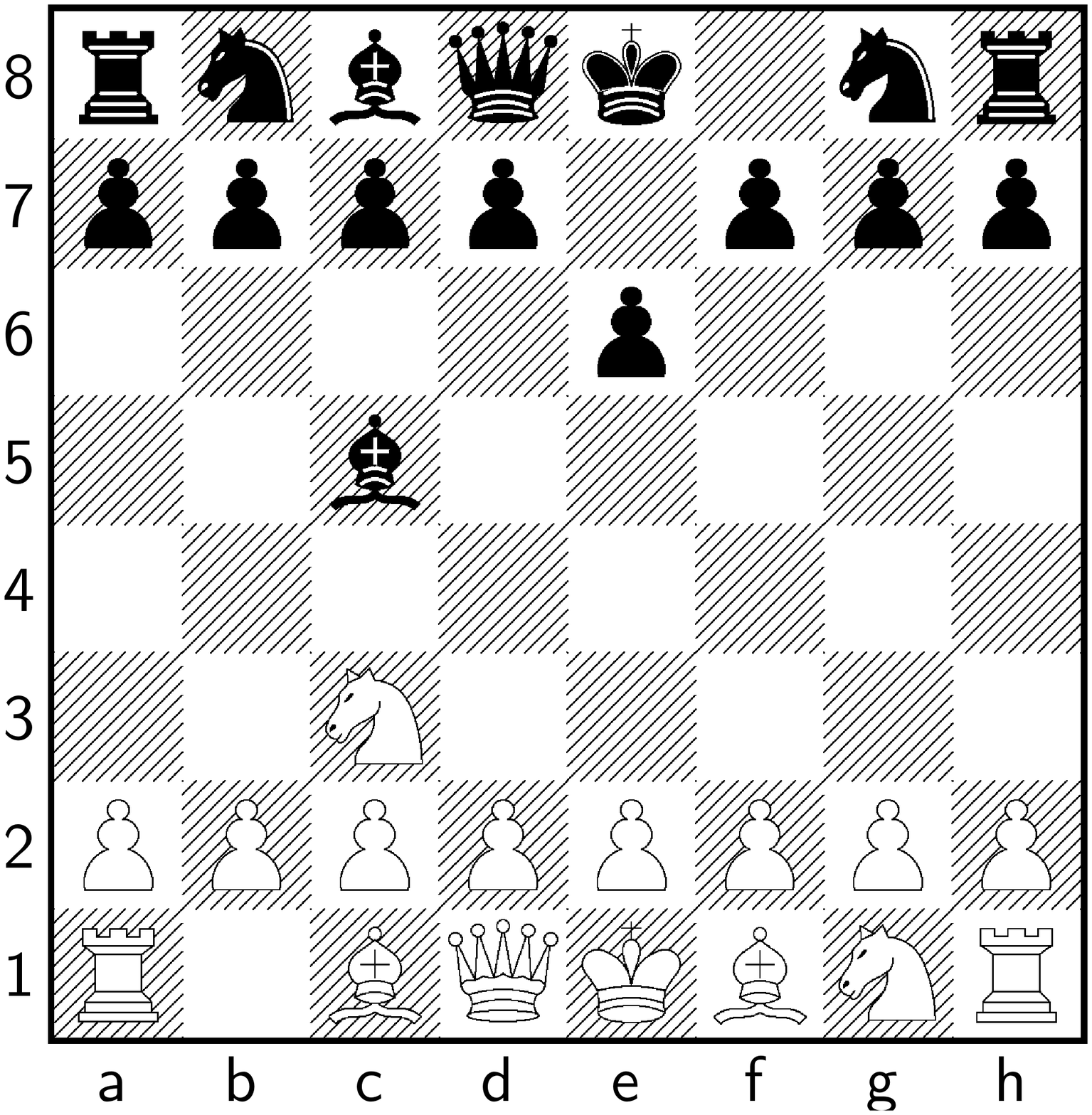}}
\end{center}
\vspace{-30pt}
Now Bob has a material advantage, but Alice has 130* chips.  Pondering the board, Bob realizes that if Alice wins the bid for less than thirty, then she can move her queen out to f6 to threaten his king and then bid everything to win the next move and take his king.  So Bob bids thirty, winning over Alice's bid of twenty-five.  He moves his knight to f3.  But Alice can still threaten Bob's king by moving her queen to h5.  Since she has 160* chips, she can now win the next two moves, regardless of what Bob bids, and capture his king.  Alice suppresses a smile as Bob realizes he has been defeated.  Head in his hands, he mumbles, \emph{``That was a total mindf**k."}

\vspace{10 pt}

{\bf Richman's Theory.}  David Richman invented and studied a class of similar bidding games in the late 1980s in which bids are allowed to be arbitrary nonnegative real numbers, not just integers.  One of Richman's discoveries is a surprising connection between such bidding games and random-turn games in which, instead of alternating moves, players flip a fair coin to determine who moves next.\footnote{The details of this result, and many other related facts, were presented by Richman's friends and collaborators in \cite{LLPU, LLPSU}}  For simplicity, suppose Alice and Bob are playing a finite, loop-free combinatorial game $G$.  Let $P(G)$ be the probability that Alice wins, assuming optimal play, and let $R(G)$ be the critical threshold between zero and one such that Alice wins if her proportion of the bidding chips is more than $R(G)$ and Bob wins if she has less than $R(G)$, with real-valued bidding.

\begin{Richman} Let $G$ be a finite combinatorial game.  Then
\[
R(G) = 1-P(G).
\]
\end{Richman}

\noindent Furthermore, a move is optimal for random-turn play if and only if it is an optimal move with real-valued bidding.  The proof of Richman's Theorem is disarmingly simple; one shows that $R(G)$ and $1-P(G)$ satisfy the same recursion with the same initial conditions, as follows.  For any position $v$ in the game $G$, let $G_v$ be the game played starting from $v$.

\begin{proof}
Suppose $v$ is an ending position, so it is a winning position for either Alice or Bob.  If $v$ is a winning position for Alice, then $R(G)$ and $1-P(G)$ are both equal to zero; if $v$ is a winning position for Bob, then both are equal to one.  Suppose $v$ is not an ending position.  By induction on the length of the game, we may assume that $R(G_w)$ is equal to $1-P(G_w)$ for every position $w$ that can be moved to from $v$.  Let $R^+(v)$ be the maximum value of $R(G_w)$, over all positions $w$ that Bob can move to from $v$, and let $R^-(v)$ be the minimum over all positions that Alice can move to.
Then it is straightforward to check that
\[
R(G_v) = \displaystyle{ \frac{R^+(v) + R^-(v)}{2}},
\]
and an optimal bid for both players is $\big| R^+(v) - R^-(v) \big| \big/2$, since Alice will always move to a position that minimizes $R$, and Bob will move to a position that maximizes it.  Similarly, Alice will always move from $v$ to a position that minimizes $1-P$, and Bob will move to a position that maximizes it.  These probabilities are equal to $R^+(v)$ and $R^-(v)$, by induction, and $1-P(G_v)$ is the average of the two, since we flip a fair coin to determine who moves next.
\end{proof}

{\bf Real vs. discrete bidding.}  Real valued bids are convenient for theoretical purposes, and are essential to Richman's elegant theory.  They are, however, disastrous for recreational play, as becomes obvious when one player bids something like $e^{-\sqrt \pi} + \log 17$.  Whole number bids are playable and fun, but the general theory and practical computation of optimal strategies become much more subtle.  For instance, the set of optimal first moves for Tic-Tac-Toe with discrete bidding depends on the number of chips in play \cite{DevelinPayne}.  Nevertheless, when the number of chips is large, discrete bidding approximates continuous bidding well enough for many purposes, and Richman's theory gives deep insight into discrete bidding game play.

{\bf Bidding Hex.} Richman's Theorem is especially exciting in light of recent developments in the theory of random-turn games.  The probabilists Peres, Schramm, Sheffield, and Wilson found an elegant solution to Random-Turn Hex, along with a Monte Carlo algorithm that quickly produces optimal or near-optimal moves \cite{PSSW07}.  This algorithm has been implemented by David Wilson \cite{Hexamania}, and the computer beats a skilled human opponent more than half the time, though anyone can beat it sometimes, by winning enough coin flips.  Elina Robeva has implemented a similar algorithm for Bidding Hex that is overwhelmingly effective---undefeated against human opponents.  

{\bf Online.}  See the Secret Blogging Seminar post by Noah Snyder \cite{Snyder08} for an excellent online introduction to bidding games, and links to further resources.  JB has developed Bidding Tic-Tac-Toe and Bidding Hex for online play through Facebook.  You can challenge your friends to a game of skill and honor, or play against the computer at a range of difficulty settings.  Visit 

\begin{center} 
\textsf{http:/\!/apps.facebook.com/biddingttt} \ and \ \textsf{http:/\!/apps.facebook.com/biddinghex}
\end{center}

\noindent and try it out!

\bibliography{math}
\bibliographystyle{amsalpha}

\end{document}